\documentclass{article}

\usepackage{amsfonts}
\usepackage[T1]{fontenc}
\usepackage{amssymb,amsthm,amsmath,mathtools,amsrefs,MnSymbol}
\usepackage{amscd,amstext,units,amsfonts,amsbsy}
\usepackage[utf8]{inputenc}
\usepackage[colorlinks=true]{hyperref}
\usepackage{color,xcolor}
\usepackage{url}
\usepackage{comment}
\usepackage{tikz-cd}
\newtheorem{theorem}{Theorem}[section]
\newtheorem{claim}{Claim}[theorem]
\newtheorem{lemma}[theorem]{Lemma}
\newtheorem{fact}[theorem]{Fact}
\newtheorem{definition}[theorem]{Definition}
\newtheorem{example}[theorem]{Example}
\newtheorem{proposition}[theorem]{Proposition}

\newtheorem{notation}[theorem]{Notation}

\newtheorem{corollary}[theorem]{Corollary}

\newcommand{\dcl}{\operatorname{dcl}}
\newcommand{\acl}{\operatorname{acl}}

\newcommand{\primes}{\operatorname{Primes}}
\newcommand{\Pres}{\operatorname{Pres}}
\newcommand{\rank}{\operatorname{rank}}

\title{Elimination of imaginaries in ordered abelian groups with bounded regular rank}
\author{Mariana Vicar\'ia}
\begin{document}
\maketitle
\begin{abstract}
In this paper we study elimination of imaginaries in some classes of pure ordered abelian groups. For the class of ordered abelian groups with bounded regular rank (equivalently with finite spines) we obtain weak elimination of imaginaries once we add sorts for the quotient groups $\Gamma/\Delta$  for each definable convex subgroup $\Delta$, and sorts for the quotient groups $\Gamma/(\Delta+ \ell\Gamma)$ where $\Delta$ is a definable convex subgroup and $\ell \in \mathbb{N}_{\geq 2}$. We refer to these sorts as the \emph{quotient sorts}. For the dp-minimal case we obtain a complete elimination of imaginaries if we also add constants to distinguish the cosets of $\ell\Gamma$ in $\Gamma$, where $\ell \in \mathbb{N}_{\geq 2}$.\\
%Finally, we study elimination of imaginaries for ordered divisible abelian groups extended with predicates to name countably many convex subgroups $(\Delta_{k}| k \in \mathbb{N})$ and obtain elimination of imaginaries once one adds the quotient groups $\Gamma/\Delta_{k}$ for each $k \in \mathbb{N}$.
\end{abstract}

\section{Introduction}
The model theory of ordered abelian groups has been studied since the sixties, and was initiated by Robinson and Zakon in \cite{RZ} who studied the completions of \emph{regular} ordered abelian groups (see Definition \ref{regular}). Later, the study of the elementary properties of ordered abelian groups was continued by Belegradek in \cite{Belegradek} for the class of \emph{poly-regular} ordered abelian groups (see Definition \ref{polyreg}). Significant achievements on (relative) quantifier elimination, model completion and definability of convex subgroups were achieved by Schmitt in \cite{supersch} for the general class of ordered abelian groups. More recently, Cluckers and Halupczok obtained a (relative) quantifier elimination for ordered abelian groups in \cite{Immi-Cluckers} in a language that is more aligned with Shelah's imaginary expansion than the one introduced by Schmitt.\\
The model theoretic classification of certain classes of ordered abelian groups is an area of active research. Results include: the well known result of Gurevich-Schmitt that no ordered abelian group has the independence property in \cite{Gurevich-Schmitt}; the dp-minimal case characterized by Jahnke, Simon and Walsberg in \cite{dpminimal}; the strongly dependent case independently obtained by Dolich-Goodrick, Farr\'e and Halevi-Hasson in \cite{Dolich-Goodrick}\cite{Farre}\cite{Halevi-Hasson} (respectively), and the distal case in  \cite{distal} due to Aschenbrenner, Chernikov, Gehret and Ziegler. \\
The next natural step regarding the model theory of ordered abelian groups was understanding a reasonable language where one will have elimination of imaginaries.
The answer to this problem, interesting in its own sake, has a significant impact in clarifying the problem of elimination of imaginaries for henselian valued fields. At the heart of the model theory of henselian valued fields is the well known \emph{Ax-Kochen/Ershov} theorem, that broadly states that the first order theory of a henselian finitely ramified valued field is completely determined by the first order theory of its residue field and its value group. In a pure henselian valued field, the value group is a pure ordered abelian group and it is interpretable in the structure.\\
Following the Ax-Kochen principle one can first attempt to solve the problem of elimination of imaginaries for henselian valued fields by following two orthogonal directions:
\begin{enumerate}
    \item The first one is to make the value group as tame as possible (e.g. to assume that it is definably complete) and to understand the obstacles that the the residue field naturally contributes to the problem. This research path was successfully finalized by Hils and Rideau-Kikuchi in \cite{AXEI}.
    \item Alternatively, one can assume the residue field is very tame (e.g. algebraically closed) and study the issues that the complexity of the value group brings to the problem. The work in \cite{multivalued} clarifies the picture for the equicharacteristic zero case, and this paper is the first milestone towards the solution. 
\end{enumerate}
To achieve elimination of imaginaries for pure ordered abelian groups, we use an abstract criterion isolated by Hrushovski in \cite{criteria} to show the following two results:
\begin{theorem}Let $\Gamma$ be an ordered abelian group of bounded regular rank (equivalently with finite spines). Then $\Gamma$ admits weak-elimination of imaginaries once the \emph{quotient sorts} are added.
\end{theorem}
\begin{theorem} Let $\Gamma$ be a dp-minimal ordered abelian group. Then $\Gamma$ admits elimination of imaginaries once the \emph{quotient sorts} are added, and we add constants to distinguish the cosets of $\ell\Gamma$ in $\Gamma$, where $\ell \in \mathbb{N}_{\geq 2}$.
\end{theorem}

%Additionally, we study a particular case of high interest in valued fields, where we consider a divisible ordered abelian group and add predicates to distinguish a countable family of convex subgroups $(\Delta_{k}\ | \ k \in \mathbb{N})$. We obtain the following result:
%\begin{thm} Let $\Gamma$ be a divisible ordered abelian group and $(\Delta_{k} \ | \ k \in \mathbb{N})$ be countable family of convex subgroups. Then $\Gamma, (\Delta_{k})_{k \in \mathbb{N}})$ admits elimination of imaginaries once one adds the quotient groups $\Gamma/\Delta_{k}$ for each $k \in \mathbb{N}$. 
%\end{thm}
%The proof for this case follows in a very similar way that the one presented for pure ordered abelian groups with finite spines, but we include it for sake of completeness as no proof is present in the existing literature.\\

 This paper is organized as follows:
 \begin{enumerate}
 \item Section \ref{preliminaries}: We present the state of model theory of ordered abelian groups and introduce the class of ordered abelian groups with bounded regular rank. 
 \item Section \ref{segments}: We characterize the definable end-segments in an ordered abelian group with bounded regular rank and show that they can be coded in the quotient sorts. 
 \item Section \ref{criterion}: We introduce Hrushovski's theorem to achieve a weak elimination of imaginaries result for the class of ordered abelian groups with bounded regular rank. This criterion requires us to check two conditions: the density of definable types, proved in Proposition \ref{density}; and
  the coding of definable types, proved in Proposition \ref{codedef}. 
 \item Section \ref{results}: We briefly present the main results for pure ordered abelian groups. 
 %\item Section \ref{OAGpredicates}: We obtain an elimination of imaginaries statement for the theory of a divisible ordered abelian groups with countably many distinguished convex subgroups. 
 \end{enumerate}

The case of direct sums of the integers with the lexicographic order has been done independently by Liccardo in \cite{Martina}, as part of her PhD thesis under D'Aquino. Hils and Mennuni in \cite{hilsmen} have independently obtained the result for the regular case. \\
\textbf{Aknowledgements:} The author would like to express her gratitude to Thomas Scanlon and Pierre Simon, for many insightful conversations and their time. The author would like to thank particularly the financial support given by the NSF Grant 1848562. The author is grateful to: R. Mennuni for his comments in a previous draft, J. Brown and A. Padgett for their careful reading, their time  and the extensive detailed grammar corrections that allowed me to  significantly improve the presentation of this paper written in my third language. Lastly, the author is also grateful to the anonymous referee for a careful reading, all the provided comments and a very fast feedback.

%\textbf{Aknowledgements:} This document was written as part of the author's  PhD thesis under  Pierre Simon and Thomas Scanlon. The author would like to express her gratitude to both of them, for many insightful conversations and their time. The author would like to thank particularly Pierre for his continuous encouragement and his financial support from his NSF Grant 1848562. The author would like to thank as well R. Mennuni for some useful comments in a previous draft. 
%The author would like as well to devote this work to those who passed away during the national protests in Colombia in May 2021,  and express her condolences to their families and friends. I thank the peaceful protestors for having taken the lead of building a more equal and diverse country, while I was working over my PhD thesis. 

\section{Preliminaries}\label{preliminaries}
\subsection*{Elimination of imaginaries}
Let $T$ be a first order theory and $\mathfrak{M}$ be its monster model. Let $D \subseteq \mathfrak{M}^{k}$ be some definable set and $E$ some definable equivalence relation over $D$. The equivalence class $e=a/E$ is said to be an \emph{imaginary element}. Imaginaries in model theory were introduced by Shelah in \cite{Shelah}. Later in \cite{Makkai}, Makkai proposed to construct the many sorted structure $\mathfrak{M}^{eq}$, where we add a sort $S_{E}$ for each definable equivalence relation $E$ and a map $\pi_{E}$ sending each element to its class. Since then, the model theoretic community has presented and studied imaginary elements in this way and refers to the multi-sorted structure $\mathfrak{M}^{eq}$ as the \emph{imaginary expansion of $\mathfrak{M}$}. We call the sorts $S_{E}$ \emph{imaginary sorts} while we refer to $\mathfrak{M}$ as the \emph{home-sort}.\\
Any formula $\phi(\mathbf{x},\mathbf{y})$ induces an equivalence relation in $\mathfrak{M}^{|\mathbf{y}|}$ defined as 
\begin{align*}
E_{\phi}(\mathbf{y}_{1}, \mathbf{y}_{2}) \ \text{if and only if}  \ \forall  \mathbf{x} \big( \phi(\mathbf{x}, \mathbf{y}_{1}) \leftrightarrow \phi(\mathbf{x},\mathbf{y}_{2}) \big).
\end{align*} 
Let $\mathbf{b} \in \mathfrak{M}^{|\mathbf{y}|}$ and $X:= \phi(\mathbf{x}, \mathbf{b})$. We call the class  $\mathbf{b}/ E_{\phi}$ the \emph{code} of $X$ and denote it as $\ulcorner X \urcorner$. We denote by $\dcl^{eq}$ and $\acl^{eq}$ the definable closure and the algebraic closure in the expansion $\mathfrak{M}^{eq}$. 
\begin{definition}
\begin{enumerate}
\item  We say that $T$ has \emph{elimination of imaginaries} if for any imaginary element $e$ there is a tuple $a$ in the home-sort such that $e \in \dcl^{eq}(a)$ and $a \in \dcl^{eq}(e)$.
\item We say that $T$ has \emph{weak elimination of imaginaries} if for any imaginary element $e$ there is a tuple $a$ in the home-sort such that $e \in \dcl^{eq}(a)$ and $a \in \acl^{eq}(e)$.
%\item We say that $T$  \emph{eliminates imaginaries}  if for any $M \vDash T$, any collection of  $M_{1},\dots,M_{k}$ of sorts in $M$, any $\emptyset$-definable set $S \subseteq M_{1}\times \dots \times M_{k}$ and any $\emptyset$-definable equivalence relation $E$ on $S$, there is an $\emptyset$-definable function $f$ from $S$ into a product of sorts, such that for any $a,b \in S$ we have:
%\begin{align*}
%E(a,b) \ \text{holds if and only if} \ f(a)=f(b).  
%\end{align*}
\item We say that $T$ \emph{codes finite sets}  if for every model $M \vDash T$ and every finite subset $S$ of $M$, the code $\ulcorner S \urcorner$ is interdefinable with a tuple of elements in $M$.
\end{enumerate}
\end{definition}
The following is a folklore fact. 
\begin{fact} \label{all} Let $T$ be a complete multi-sorted theory.  If $T$ has weak elimination of imaginaries and codes finite sets then $T$ eliminates imaginaries.
\end{fact}
\textbf{Ordered abelian groups of bounded regular rank}\label{background}\\
In this section we summarize several results about the classification of ordered abelian groups and their model theoretic behavior. We start by recalling the following folklore fact. 
\begin{fact} Let $(\Gamma,< ,+, 0)$ be a non-trivial ordered abelian group. Then the topology induced by the order in $\Gamma$ is discrete if and only if $\Gamma$  has a minimal positive element. In this case we say that $\Gamma$ is \emph{discrete}, otherwise we say that it is \emph{dense}. 
\end{fact}

The following notions were isolated in the sixties by Robinson and Zakon in \cite{RZ} to understand some complete extensions of the theory of ordered abelian groups.
\begin{definition} Let $\Gamma$ be an ordered abelian group and $n \in \mathbb{N}_{\geq 2}$.
\begin{enumerate}
\item  Let $\gamma \in \Gamma$. We say that $\gamma$ is \emph{$n$-divisible} if there is some $\beta \in \Gamma$ such that $\gamma=n\beta$. 
\item We say that $\Gamma$ is \emph{$n$-divisible} if every element $\gamma \in \Gamma$ is $n$-divisible. 
\item $\Gamma$ is said to be \emph{$n$-regular} if any interval with at least $n$ points contains an $n$-divisible element. 
\end{enumerate}
\end{definition}
\begin{example} We include some examples to illustrate the previous definitions.
\begin{enumerate}
    \item Consider the ordered abelian group $(\mathbb{Z},+, \leq, 0)$, the elements $2,4,6$ are $2$-divisible while $1$ is not. 
    \item The groups $(\mathbb{Q},+,\leq, 0)$ and $(\mathbb{Z}, +, \leq, 0)$ are $n$-regular for each natural number $n \in \mathbb{N}_{\geq 2}$. The group  $(\mathbb{Z} \oplus \mathbb{Z}, \leq_{lex}, +, 0)$, where $\leq_{lex}$ is the lexicographic order, is not $2$-regular, because the interval $\big((1,-1),(1,4)\big)=\{ (1,0), (1,1), (1,2),(1,3)\}$ does not contain a $2$-divisible element.
\end{enumerate}
\end{example}

The following definitions were introduced by Schmitt in \cite{sch} and \cite{supersch}.
\begin{definition}  We fix an ordered abelian group $\Gamma$ and $n\in \mathbb{N}_{\geq 2}$. Let $\gamma \in \Gamma$. We define:
\begin{itemize}
    \item $A(\gamma)=$ the largest convex subgroup of $\Gamma$ not containing $\gamma$.
    \item $B(\gamma)=$ the smallest convex subgroup of $\Gamma$ containing $\gamma$.
    \item $C(\gamma)=B(\gamma)/A(\gamma)$.
    \item $A_{n}(\gamma)=$ the smallest convex subgroup $C$ of $\Gamma$ such that $B(g)/C$ is $n$-regular. 
    \item $B_{n}(g)=$ the largest convex subgroup $C$ of $\Gamma$ such that $C/A_{n}(\gamma)$ is $n$-regular.
\end{itemize}
\end{definition}
In \cite[Chapter 2]{supersch}, Schmitt shows that the groups $A_{n}(\gamma)$ and $B_{n}(\gamma)$ are definable in the language of ordered abelian groups $\mathcal{L}_{OAG}=\{ +,-, \leq, 0\}$ by a first order formula using only the parameter $\gamma$.\\
We recall that the set of convex subgroups of an ordered abelian group is totally ordered by inclusion.
\begin{definition}
Let $\Gamma$ be an ordered abelian group and $n \in \mathbb{N}_{\geq 2}$, we define the \emph{$n$-regular rank} to be the order type of: \\
\begin{align*}
\big( \{ A_{n}(\gamma) \ | \ \gamma \in \Gamma \backslash \{ 0\}\}, \subseteq \big).
    \end{align*}
\end{definition}
The $n$-regular rank of an ordered abelian group $\Gamma$ is a linear order, and when it is finite we can identify it with its cardinal. In \cite{Farre}, Farr\'e emphasizes that we can characterize it without mentioning the subgroups $A_{n}(\gamma)$. The following is \cite[Remark 2.2]{Farre}. 
\begin{definition}
Let $\Gamma$ be an ordered abelian group and $n \in \mathbb{N}_{\geq 2}$, then: 
\begin{enumerate}
\item $\Gamma$ has $n$-regular rank equal to $0$ if and only if $\Gamma=\{0\}$,
\item $\Gamma$ has $n$-regular rank equal to $1$ if and only if $\Gamma$ is $n$-regular and not trivial, 
\item $\Gamma$ has $n$-regular rank equal to $m$ if there are $\Delta_{0}, \dots, \Delta_{m}$ convex subgroups of $\Gamma$, such that:
\begin{itemize}
\item $\{0\}=\Delta_{0} < \Delta_{1} < \dots < \Delta_{m}=\Gamma$,
\item for each $0 \leq i< m$, the quotient group $\Delta_{i+1}/\Delta_{i}$ is $n$-regular,
 \item the quotient group $\Delta_{i+1}/\Delta_{i}$ is not $n$-divisible for $0 < i < m$. 
\end{itemize}
In this case we define $RJ_{n}(\Gamma)=\{ \Delta_{0}, \dots, \Delta_{m-1}\}$. The elements of this set are called the \emph{$n$-regular jumps}.
\end{enumerate}
\end{definition}
\begin{example}
Let $G=\underbrace{\mathbb{Z} \oplus \dots \oplus \mathbb{Z}}_{n- \text{times}}$ with the lexicographic order $\leq_{lex}$. The $3$-regular rank of $G$ is equal to $n$. This is witnessed by the sequence:
\begin{align*}
\{\mathbf{0}\} \leq \underbrace{\{0\} \oplus \dots \oplus \{0\}}_{n-1\ \text{times}} \oplus \mathbb{Z}\leq\dots \leq 
\{0\}\oplus \underbrace{\mathbb{Z} \oplus \dots \oplus \mathbb{Z}}_{n-1- \text{times}} \leq  \mathbb{Z} \oplus \dots \oplus \mathbb{Z}.
\end{align*}
\end{example}

\subsubsection*{Regular groups and poly-regular groups}
\begin{definition}\label{regular} An ordered abelian group $\Gamma$ is said to be \emph{regular} if it is $n$-regular for all $n \in \mathbb{N}$.
\end{definition}
\begin{example}
$\displaystyle{(\mathbb{Z}, +,\leq, 0)}$ and $\displaystyle{(\mathbb{Q},+,\leq, 0)}$  are standard examples of regular groups. By \cite[Theorem 1.2]{Belegradek} any archimedean group is regular. 
\end{example}
Robinson and Zakon in their seminal paper \cite{RZ} completely characterized the possible completions of the theory of regular groups, obtained by extending the first order theory of ordered abelian groups with axioms asserting that for each $n \in \mathbb{N}$ 
if an interval contains at least $n$-elements then it contains an $n$-divisible element. 
The following is \cite[Theorem 4.7]{RZ}. 
\begin{theorem}\label{charreg}
The possible completions of the theory of regular groups, are: 
\begin{enumerate}
    \item the theory of discrete regular groups, and
    \item  the completions of the theory of dense regular groups $T_{\chi}$ where
    \begin{align*}
      \chi&=:\primes \rightarrow \mathbb{N} \cup \{\infty\},\\
    \end{align*}
   is a function specifying the index $\chi(p)= [\Gamma : p\Gamma]$.
\end{enumerate}
\end{theorem}
Robinson and Zakon proved as well that each of these completions is the theory of some archimedean group. In particular, any discrete regular group is elementarily equivalent to $(\mathbb{Z},\leq, +, 0)$. This theory is called the theory of Presburger arithmetic, introduced in $1929$ by M. Presburger, who proved that it admits quantifier elimination in the well known \emph{Presburger Language}  $\mathcal{L}_{\Pres}=\{ 0,1,+, -, <, (\equiv_{m})_{m \in \mathbb{N}_{\geq 2}}\}$. Given an ordered abelian group $\Gamma$ we naturally see it as a $\mathcal{L}_{\Pres}$-structure. The symbols $\{0,+,- ,<\}$ take their obvious interpretation. If $\Gamma$ is discrete, the constant symbol $1$ is interpreted as the least positive element of $\Gamma$,  and by $0$ otherwise. For each $m \in \mathbb{N}_{\geq 2}$ the binary relation symbol $\equiv_{m}$ is interpreted as the equivalence modulo $m$, i.e. for any $g,h \in \Gamma$ $\displaystyle{g \equiv_{m} h \ \text{if and only if \ } g-h \in m\Gamma}$. \\
 The theory of a dense ordered abelian group admits quantifier elimination in the Presburger language if and only if it is regular. This is a result of Weispfenning in \cite[Theorem 2.9]{poly}.
\begin{definition}\label{polyreg} Let $\Gamma$ be an ordered abelian group. We say that it is \emph{poly-regular} if it is elementarily equivalent to a subgroup of the lexicographically ordered group $\mathbb{R}^{n}$.
\end{definition}
In \cite{Belegradek} Belegradek studied poly-regular groups and proved that an ordered abelian group is poly-regular if and only if it has finitely many proper definable convex subgroups, and all the proper definable subgroups are definable over the empty set. In \cite[Theorem 2.9]{poly} Weispfenning obtained quantifier elimination for the class of poly-regular groups in the language of ordered abelian groups extended with predicates to distinguish the subgroups $\Delta+\ell \Gamma$ where $\Delta$ is a convex subgroup and $\ell \in \mathbb{N}_{\geq 2}$. 
\subsubsection*{Ordered abelian groups with bounded regular rank}
\begin{definition}
Let $\Gamma$ be an ordered abelian group. We say that it has \emph{bounded regular rank} if it has finite $n$-regular rank for each $n \in \mathbb{N}_{\geq 2}$.  For notation, we will use $\displaystyle{RJ(\Gamma)= \bigcup_{n \in \mathbb{N}_{\geq 2}} RJ_{n}(\Gamma)}$. 
\end{definition}

The class of ordered abelian groups of bounded regular rank  extends the class of poly-regular groups and regular groups. The terminology of \emph{bounded regular rank} becomes clear with the following Proposition (item $3$). 
\begin{proposition}\label{charbounded} Let $\Gamma$ be an ordered abelian group. The following are all equivalent: 
\begin{enumerate}
\item $\Gamma$ has finite $p$-regular rank for each prime number $p$.
\item $\Gamma$ has finite $n$-regular rank for each $n \geq 2$.
\item There is some cardinal $\kappa$ such that for any $H \equiv \Gamma$, $|RJ(H)|\leq \kappa$.
\item For any $H \equiv \Gamma$, any definable convex subgroup of $H$  has a definition without parameters.
\item There is some cardinal $\kappa$ such that for any $H \equiv \Gamma$, $H$ has at most $\kappa$  definable convex subgroups. 
\end{enumerate}
Moreover, in this case $RJ(\Gamma)$ is the collection of all proper definable convex subgroups of $\Gamma$ and all are definable without parameters. In particular, there are only countably many definable convex subgroups. 
\end{proposition}
\begin{proof}
This is \cite[Proposition 2.3]{Farre}. 
\end{proof}
\subsubsection*{Quantifier elimination and the quotient sorts}
 In \cite{Immi-Cluckers} Cluckers and Halupczok introduced a language $\mathcal{L}_{qe}$  to obtain  quantifier elimination for ordered abelian groups relative to the \emph{auxiliary sorts} $S_{n}$, $T_{n}$ and $T_{n}^{+}$, whose precise description can be found in \cite[Definition 1.5]{Immi-Cluckers}. This language is similar in spirit to the one introduced by Schmitt in  \cite{supersch}, but has lately been preferred by the community as it is more in line with the many-sorted language of Shelah's imaginary expansion $\mathfrak{M}^{eq}$. Schmitt does not distinguish between the sorts $S_{n}$, $T_{n}$ and $T_{n}^{+}$. Instead for each $n \in \mathbb{N}$ he works with a single sort $Sp_{n}(\Gamma) $ called the \emph{$n$-spine} of $\Gamma$, whose description can be found in  \cite[Section 2]{Gurevich-Schmitt}. In  \cite[Section 1.5]{Immi-Cluckers} it is explained how the auxiliary sorts of Cluckers and Halupczok are related to the $n$-spines $Sp_{n}(\Gamma)$ of Schmitt. In \cite[Section 2]{Farre}, it is shown that an ordered abelian group $\Gamma$ has bounded regular rank if and only if all the $n$-spines are finite, and $Sp_{n}(\Gamma)=RJ_{n}(\Gamma)$. In this case, we define the regular rank of $\Gamma$ as the cardinal $|RJ(\Gamma)|$, which is either finite or $\aleph_{0}$. Instead of saying that $\Gamma$ is an ordered abelian group with finite spines, we prefer to use the classical terminology of bounded regular rank, as it emphasizes the relevance of the $n$-regular jumps and the role of the divisibilities to describe the definable convex subgroups.

\begin{definition} [The language $\mathcal{L}$]
Let $\Gamma$ be an ordered abelian group with bounded regular rank. We view $\Gamma$ as a multi-sorted structure in the language $\mathcal{L}$, where: 
\begin{enumerate}
\item we add a sort for the ordered abelian group $\Gamma$, and we equip it with a copy of the language $\mathcal{L}_{\Pres}$ extended with predicates to distinguish each of the convex subgroups $\Delta \in RJ(\Gamma)$. We refer to this sort as the \emph{main sort}. 
\item For each $\Delta \in RJ(\Gamma)$ we add a sort for each of the ordered abelian groups $\Gamma/\Delta$, equipped with a copy of the Presburger language\\ 
$\mathcal{L}_{\Pres}^{\Delta}=\{ 0^{\Delta}, 1^{\Delta}, +^{\Delta}, -^{\Delta}, <^{\Delta}, (\equiv_{m}^{\Delta})_{m \in \mathbb{N}_{\geq 2}}\}$. \\
We add as well a map $\rho_{\Delta}: \Gamma \rightarrow \Gamma/\Delta$, interpreted as the natural quotient map. 
\end{enumerate}
\end{definition}
In \cite[Theorem 2.4]{Farre} Farr\'e obtained a quantifier elimination statement for the class of ordered abelian groups with bounded regular rank in the languange $\mathcal{L}$ extended with a set of constants in the home sort. However, we present a slightly different language where we add the constants for the minimal element in $\Gamma/\Delta$ (if it exists) instead of adding a representative in the home-sort whose projection is the minimal class in $\Gamma/\Delta$.  For this purpose we highlight that the following statement is a direct consequence of \cite[Proposition 3.14]{distal}. 
\begin{theorem} \label{QEbounded} Let $\Gamma$ be an ordered abelian group with bounded regular rank. Then $\Gamma$ admits quantifier elimination in the language $\mathcal{L}$. 
\end{theorem}
\begin{notation}
We will be mainly interested in the description of the definable sets in the \emph{main sort}. For this purpose we will slightly abuse the language, to simplify the notation. For each $k \in \mathbb{Z}$ and $\Delta \in RJ(\Gamma)$ we define $k^{\Delta}= k \cdot 1^{\Delta}$, where $1^{\Delta}$ is the minimal element in $\Gamma/\Delta$ if it exists. We will sometimes indicate $k^{\Delta}$ simply as $k+\Delta$. We introduce the following notation:
\begin{enumerate}
\item We write $\tau(\mathbf{x}) + \Delta < \beta + k+ \Delta$ for the formula $\rho_{\Delta}(\tau(\mathbf{x}))<^{\Delta} \rho_{\Delta}(\beta)+k^{\Delta}$.
\item We write $\tau(\mathbf{x}) \equiv_{\Delta} \beta +k$ for the formula $\rho_{\Delta}(\tau(\mathbf{x}))= \rho_{\Delta}(\beta)+k^{\Delta}$.
\item We write $\tau(\mathbf{x}) \equiv_{\Delta + m\Gamma} \beta+ k$ for the formula $\rho_{\Delta}(\tau(\mathbf{x})) \equiv_{m}^{\Delta} \rho_{\Delta}(\beta)+ k^{\Delta}$. The latter is interpreted as  $\rho_{\Delta}(\tau(\mathbf{x})) - (\rho_{\Delta}(\beta)+k^{\Delta}) \in m \big(\Gamma/ \Delta)$. 
\end{enumerate}
Here$\tau(\mathbf{\mathbf{x}})$ is a term in the language of ordered abelian groups in $m$ variables, $\mathbf{x}=(x_{1}, \dots, x_{m})$ and  $\beta \in \Gamma$. 
\end{notation}
\begin{definition}
\begin{enumerate}
\item A set $S \subset \Gamma$ is said to be an \emph{end-segment} (respectively an initial segment) if for any $x \in S$ and $y \in \Gamma$, $x < y$ (respectively $y < x$) we have that $y \in S$. 
\item Let $n \in \mathbb{N}_{\geq 2}$, $\Delta \in RJ(\Gamma)$, $\beta \in \Gamma \cup \{ -\infty, +\infty\}$ and $\square \in \{ \geq, >\} $.\\ 
$\{ \eta \in \Gamma \ | \ n\eta + \Delta \square  \beta +\Delta\}$ 
is an end-segment of $\Gamma$. We call the end-segments of this form \emph{divisibility end-segments}. We define \emph{divisibility initial segments} analogously.
 
\item A \emph{mid-segment} is a non empty set $C$ of the form $C= U \cap L$ where $U$ is a divisibility end-segment and $L$ is a divisibility initial segment.
\item A \emph{basic positive congruence formula} is a formula of the form $zx \equiv_{\Delta+l\Gamma} \beta + k$ where $\beta \in \Gamma$, $z,k \in \mathbb{Z}$ and $l \in \mathbb{N}_{\geq 2}$. Likewise,  a \emph{basic negative formula} is a formula of the form  $zx \not \equiv_{\Delta+l\Gamma} \beta + k$. A \emph{basic congruence formula} is either a basic positive congruence formula or a basic negative formula. 
\item A \emph{finite congruence restriction} is a finite conjunction of basic congruence formulas. 
\item A \emph{nice set} is a set of the form $C \cap X$, where $C$ is a mid-segment and $X$ is defined by a finite congruence restriction. 
\end{enumerate}
\end{definition}
The following is a direct consequence of quantifier elimination.
\begin{corollary}\label{nice} Let $\Gamma$ be an ordered abelian group with bounded regular rank. Let $X \subseteq \Gamma$ be a definable set. Then $X$ is a finite union of nice sets. 
\end{corollary}
We will consider an extension $\mathcal{L}_{Q}$ of the language $\mathcal{L}$, where for each natural number $n \in \mathbb{N}_{\geq 2}$ and $\Delta \in RJ(\Gamma)$ we add a sort for the quotient group $\Gamma/ (\Delta+n\Gamma)$ and a map $\pi_{\Delta}^{n}:\Gamma \rightarrow \Gamma/(\Delta+n\Gamma)$. 

We will refer to the sorts in the language $\mathcal{L}_{Q}$ as \emph{quotient sorts}. \\

The following fact will be very useful to show weak elimination of imaginaries for ordered abelian groups with bounded regular rank.
\begin{fact}\label{quotientsareok} Let $\Gamma$ be an ordered abelian group of finite  $n$-regular rank witnessed by the sequence $\{ 0 \}=\Delta_{0} < \Delta_{1} < ...< \Delta_{l}=\Gamma$ and fix some definable convex subgroup  $H$. Then $\Gamma/H$ is also a group of finite $n$-regular rank. Moreover, if $\Gamma$ is an ordered abelian group of bounded regular rank, then $H \in RJ(\Gamma)$ and each coset of  $\Delta_{i}/H$ in $\Gamma/H$ is interdefinable with an element of $\Gamma/\Delta_{i}$. 
\end{fact}
\begin{proof}
Let $\Gamma$ be an ordered abelian group and $H$ a convex subgroup. Assume that $\Gamma$ has finite $n$-regular rank, witnessed by the sequence $\{ 0 \}=\Delta_{0} < \Delta_{1} < ...< \Delta_{l}=\Gamma$ and let $r$ be the smallest index such that $\Delta_{r} \subseteq H \subseteq \Delta_{r+1}$. We aim to show that $\Delta_{r+1}/H < \dots < \Delta_{l}/H=\Gamma/H$ witnesses that $\Gamma/H$ has finite $n$-regular rank. For each $r\leq i< l$, by the isomorphism theorem $(\Delta_{i+1}/H) / (\Delta_{i}/H) \cong \Delta_{i+1}/\Delta_{i}$. As $\Delta_{i+1}/\Delta_{i}$ is $n$-regular and not $n$-divisible, so is $(\Delta_{i+1}/H) / (\Delta_{i}/H)$.\\
The second part of the statement follows immediately by the isomorphism theorem. $\square$
\end{proof}

\subsection*{A survey of model theoretic results on ordered abelian groups}
In $1984$ the classification of the model theoretic complexity of ordered abelian groups was initiated by Gurevich and Schmitt in \cite{Gurevich-Schmitt}, who proved that no ordered abelian group has the independence property. During the last years finer classifications have been achieved, and we present the state of the field in this subsection. 
\begin{definition} Let $\Gamma$ be an ordered abelian group and let $p$ be a prime number. We say that $p$ is a \emph{singular prime} if $[\Gamma: p\Gamma]= \infty$. If $\Gamma$ does not have singular primes we call it \emph{non-singular}. 
\end{definition}
The following result corresponds to \cite[Proposition 5.1]{dpminimal}.
\begin{proposition}\label{dpmin} Let $\Gamma$ be an ordered abelian group. The following conditions are equivalent:
\begin{enumerate}
\item $\Gamma$ is non-singular,
\item $\Gamma$ is $dp$-minimal.
\end{enumerate}
\end{proposition}
The following is \cite[ Theorem 3.13]{distal}.
\begin{proposition} Let $\Gamma$ be an ordered abelian group with bounded regular rank (i.e. each $Sp_{n}(\Gamma)$ is finite). The following statements are equivalent:
\begin{enumerate}
\item $\Gamma$ is distal,
\item $\Gamma$ is $dp$-minimal. 
\end{enumerate}
\end{proposition}
The following statement was independently achieved in  \cite{Dolich-Goodrick},\cite{Farre} and \cite{Halevi-Hasson}. 
\begin{proposition} Let $\Gamma$ be an ordered abelian group. The following conditions are equivalent:
\begin{enumerate}
\item $\Gamma$ is strongly dependent,
\item $\Gamma$ has finite $dp$-rank,
\item $\Gamma$ has bounded regular rank and finitely many singular primes. 
\end{enumerate}
Moreover, let $\mathcal{P}= \{ p \in \mathbb{N} \ | \ p$ is a singular prime$\}$. Then
\begin{align*}
dp-\rank(\Gamma) \leq 1+ \sum_{p \in \mathcal{P}} |RJ_{p}(G)|. 
\end{align*}
\end{proposition}
\section{Definable end-segments}\label{segments}
In this subsection we characterize the definable end-segments (or initial segments) in an ordered abelian group with bounded regular rank (equivalently with finite spines). We also show that they can be coded in the quotient sorts. 
\begin{definition} Let $\Gamma$ be an ordered abelian group with bounded regular rank: 
\begin{enumerate}
\item Given $S \subseteq \Gamma$ an end-segment (or an initial segment) we denote by $\Delta_{S}$ the stabilizer of $S$, i.e. $\Delta_{S}:= \{ \beta \in \Gamma \ | \ \beta+S= S\}$.
\item Let $S \subseteq \Gamma$ be an end-segment and $\Delta \in RJ(\Gamma)$. We consider the projection map $\displaystyle{\rho_{\Delta}: \Gamma \rightarrow \Gamma/\Delta}$, and we will denote $\rho_{\Delta}(S)$
by $S_{\Delta}$. One can show that $\displaystyle{S_{\Delta}=\{ \gamma \in \Gamma/\Delta \ | \ \exists y \in S \  \rho_{\Delta}(y)=\gamma\}}$ is a definable end-segment of $\Gamma/\Delta$, as it is the projection of an end-segment. 
\item Let $\Delta \in RJ(\Gamma)$ and $S \subseteq \Gamma$ be an end-segment. 
We say that $S$ is $\Delta$-\emph{decomposable} if it is a union of $\Delta$-cosets. 
\item Let $X$ and $Y$ be definable sets. We say that $Y$ is \emph{coinitial} (or \emph{cofinal}) in $X$ if for any $y \in X$
 there is some element $z \in X \cap Y$ such that $z \leq  y$ (respectively $z \geq y$).
\end{enumerate}
\end{definition}
\begin{fact}\label{stab}Let $S \subseteq \Gamma$ be a definable end-segment. Then $\Delta_{S}$ is
 a definable convex subgroup of $\Gamma$, and therefore $\Delta_{S} \in RJ(\Gamma)$. Furthermore, 
 $\displaystyle{\Delta_{S}= \bigcup_{\Delta \in \mathcal{C}} \Delta}$, where 
 \begin{align*}
 \mathcal{C}=\{\Delta \in RJ(\Gamma) \ | \ S \text{\ is $\Delta$-decomposable}\}. 
 \end{align*}
\end{fact}
\begin{proof}
We first show that $\displaystyle{\Delta_{S} \subseteq \bigcup_{\Delta \in \mathcal{C}} \Delta}$. Note $\Delta_{S}$ is a definable convex subgroup, so $\Delta_{S} \in RJ(\Gamma)$. We aim to show that $S$ is $\Delta_{S}$-decomposable, so it is sufficient to show that for any $\gamma \in S$, $\gamma+\Delta_{S} \subseteq S$. Fix some $\gamma \in \Gamma$. If $\delta \in \Delta_{S}$ then $\gamma+\delta \in S$, so $\gamma+ \Delta_{S}\subseteq S$.\\
We now prove that $\displaystyle{\bigcup_{\Delta \in \mathcal{C}}} \Delta \subseteq \Delta_{S}$. Let $\Delta \in \mathcal{C}$ and fix some $\delta \in \Delta$. We want to show that $\delta+S \subseteq S$ and $S \subseteq \delta+S$. 
Because $S$ is $\Delta$-decomposable, $\gamma+\Delta \subseteq S$ for any $\gamma \in S$. In particular $\gamma+\delta \in S$. As $\gamma$ is an arbitrary element in $S$, we conclude that $\delta+S \subseteq S$. It is only left to show that $S \subseteq \delta +S$. Let $\gamma \in S$, then $\gamma-\delta \in S$ because $\gamma+\Delta \subseteq S$. As $\gamma=\delta+(\gamma-\delta)\in \delta +S$, we have $S \subseteq \delta+S$, as required. $\square$ 
\end{proof}
\begin{proposition}\label{charend} Let $\Gamma$ be an ordered abelian group of bounded regular rank. 
Any definable end-segment is a divisibility end-segment. 
\end{proposition}
\begin{proof}
Let $S \subseteq \Gamma$ be a definable end-segment such that $S \neq \Gamma$. By Fact \ref{stab}, $\Delta_{S}$ is a definable convex subgroup of $\Gamma$ and  $S$ is $\Delta_{S}$-decomposable. 
To simplify the notation we will denote $\hat{\Gamma}=\Gamma/\Delta_{S}$ and $\hat{S}=S_{\Delta_{S}}= \rho_{\Delta_{S}}(S)$. It is sufficient to prove that $\hat{S}$ is a divisibility end-segment in $\hat{\Gamma}$. \\
\begin{claim}
Note that for any $k \in \mathbb{N}$ exactly one of the following  occurs: 
\begin{itemize}
\item $\hat{\Gamma}$ is $k$-regular.
\item There is a non trivial $k$-regular convex subgroup $\Lambda_{k}$ of $\hat{\Gamma}$ and a coset $\eta+\Lambda_{k}$ such that $\hat{S} \cap (\eta+ \Lambda_{k}) \neq \emptyset$ and $\hat{S}^{c} \cap (\eta+ \Lambda_{k}) \neq \emptyset$.
\end{itemize}
\end{claim}
\begin{proof}
Let  $\{0\}=\Delta_{0} < \Delta_{1} < \dots < \Delta_{l}=\Gamma$ the sequence of convex subgroups witnessing that $\Gamma$ has $k$-finite regular rank equal to $l$. Let $r \leq l$ be the smallest index such that $\Delta_{S} \subsetneq \Delta_{r}$. If $r=l$ then $\hat{\Gamma}$ is $k$-regular. Otherwise the quotient group $\Lambda_{k}=\Delta_{r}/ \Delta_{S}$ satisfies the required conditions. Indeed,  as $\Delta_{r}/\Delta_{r-1}$ is $k$-regular so is $\Delta_{r}/ \Delta_{S}$. Additionally, $S$ is not $\Delta_{r}$-decomposable. If it were, then we would have $\Delta_{r} \subseteq \Delta_{S}$ which contradicts $\Delta_{S} \subsetneq \Delta_{r}$.  Then there is some coset $\eta+ \Delta_{r}$ such that $S \cap (\eta+\Delta_{r}) \neq \emptyset$ and $S^{c} \cap (\eta+\Delta_{r}) \neq \emptyset$ because otherwise $S$ would be $\Delta_{r}$-decomposable. Thus $\hat{S} \cap (\hat{\eta}+\Lambda_{k} )\neq \emptyset$ and $\hat{S}^{c} \cap (\hat{\eta}+ \Lambda_{k} )\neq \emptyset$, where $\hat{\eta}= \eta+ \Delta_{S}$.\  $\square$
\end{proof}

We may assume that $\hat{S}$ does not have a minimum because otherwise the statement follows immediately. By Corollary \ref{nice} applied to $\hat{\Gamma}$, $\hat{S}$ is a finite union of nice sets $C_{i} \cap X_{i}$, where $C_{i}=U_{i}\cap L_{i}$. 
%\begin{enumerate}
%\item Each $C_{i}$ is a mid-segment defined as the intersection $U_{i} \cap L_{i}$, where $U_{i}$ is a divisibility end-segment and $L_{i}$ is a divisibility initial segment. 
%\item Each $X_{i}$ is a finite restriction of congruences i.e. 
%\begin{align*}
%X_{i}(x):= \bigwedge_{j \leq k_{i}} \big(n_{j}x \equiv_{\hat{\Delta}_{j}+ m_{j}\hat{\Gamma}} \tau_{j}(\beta_{j}) \big)  \wedge \bigwedge_{l \leq w_{i}} {k_{l}x \nequiv_{\hat{\Delta}_{l}+ s_{l}\hat{\Gamma}} \tau_{j}(\eta_{l}) \big)}.
%\end{align*}
%\end{enumerate}
%where  each $n_{j}, k_{l} \in \mathbb{Z}$, $\hat{\Delta}_{j},\hat{\Delta}_{l}  \in RJ(\hat{\Gamma})$, $m_{j}, s_{l} \in \mathbb{N}^{>0}$,  $\tau_{j}(\mathbf{y})$ are terms in the language of ordered abelian groups (increased with possible constants) and $\beta_{j},\beta_{l}$ are tuples of elements in $\hat{\Gamma}$. \\
As $\hat{S}$ is a definable end-segment, it is sufficient to understand the co-initial description of $\hat{S}$. Without loss of generality we may assume that $U_{i} \subseteq U_{1}$ for all $i$. Let $\hat{\Delta} \in RJ(\hat{\Gamma})$. Then $\hat{\Delta}=\Delta/\Delta_{S}$ for some $\Delta_{S} \subsetneq \Delta \in RJ(\Gamma)$. Thus there is a coset $\eta+ \hat{\Delta}$ such that 
 $\hat{S} \cap (\eta+ \hat{\Delta}) \neq \emptyset$ and $\hat{S}^{c} \cap (\eta+ \hat{\Delta}) \neq \emptyset$ because $S$ is not $\Delta$-decomposable. 
 
 Hence, each of the congruence formulas involving the groups $\hat{\Delta}+ k \hat{\Gamma}$ does not change its truth value over $U_{1}  \cap (\eta+ \hat{\Delta})$. Therefore it does not change its truth value co-initially in $U_{1}$. \\
Consider a conjunction of congruence restrictions of the form: 
\begin{align*}
C(x):=\big( \bigwedge_{i\leq s} x \equiv_{k_{i}\hat{\Gamma} } c_{i} \big)\wedge \big( \bigwedge_{j \leq l }  \neg(x \equiv_{r_{j} \hat{\Gamma}} d_{j}) \big).
\end{align*}
Let $M$ be the least common multiple of all the $k_{i}$'s and $r_{j}$'s involved in the definition of $C(x)$. By the previous Claim, $\Gamma$ is $M$-regular or we can find an $M$-regular group $\Lambda_{M}$ and a coset that intersects $\hat{S}$ and its complement. We first assume the existence of a non-trivial convex subgroup $\Lambda_{M}$ and a coset $\eta+\Lambda_{M}$  such that $\hat{S} \cap (\eta+ \Lambda_{M}) \neq \emptyset$ and $\hat{S}^{c} \cap (\eta+ \Lambda_{M}) \neq \emptyset$. Let $Y=C(x) \cap (U_{1} \cap (\eta+ \Lambda_{M}) )$.\\

\begin{claim}{If $Y \neq \emptyset$, then $C(x)$ is co-initial in $U_{1}$.}\end{claim}
\begin{proof}
Let $x_{0} \in Y$ and $ U'= \big(U_{1} \cap (\eta+ \Lambda_{M})\big)-x_{0}$. $U'$ is a definable end-segment of $\Lambda_{M}$ without a minimum. Fix an element $\delta \in U'$. As $U'$ does not have a minimum and $\Lambda_{M}$ is $M$-regular, we can find an element $\gamma \in \Lambda_{M}$ such that $M\gamma \in U'$ and $M\gamma<\delta$. Then $z= M\gamma+x_{0} \in Y$ and $z< x_{0}+\delta$. Thus $C(x)$ is co-initial in $U_{1}$.\ $\square$
\end{proof}
Likewise, if $\hat{\Gamma}$ is $M$-regular we can conclude that $C(x)$ is co-initial in $U_{1}$. Consequently, the congruence restrictions are irrelevant in the definition of the end-segment $S$. It must be the case then that $S=U_{1}$, as desired.\ $\square$ 
\end{proof}
Though that it may seem like any divisibility cut defined by a formula of the form $nx \square \beta$  where $n \in \mathbb{N}_{\geq 2}$, $\square \in \{ \geq, >\}$ and $\beta \in \Gamma$  could be coded by $\beta$, this statement is false and requires a slightly more delicate treatment. We introduce the following example to motivate the reader to not dismiss the technical work in Proposition \ref{charend}. 
\begin{example} Consider the ordered abelian group $(\mathbb{Z}\oplus \mathbb{Z},\leq_{lex},+,0)$ where $\leq_{lex}$ is the lexicographic order. The definable end-segment $\displaystyle{S=\{ z \in \mathbb{Z}^{2} \ | \ 2z \geq (1,1)\}}$. Note that for any $\beta \in \mathbb{Z}$, $S$ is also defined by the formula $2z \geq (1, \beta)$. 
\end{example}

\begin{lemma}\label{cociente} Let $\Gamma$ be an ordered abelian group of bounded regular rank.  Let $\{0\}=\Delta_{0} \leq \Delta_{1} \leq \dots \leq \Delta_{l}=\Gamma$ be the sequence of convex subgroups witnessing that $\Gamma$ has finite $n$-regular rank. Then any divisibility end-segment $S$ defined by a formula $nx \ \square \ \beta$ where $n \in \mathbb{N}_{\geq 1}$, $\square \in \{ \geq, >\}$ and $\beta \in \Gamma$  is coded by a tuple of elements in the sorts $\Gamma \cup \{ \Gamma/\Delta_{i} \ | \ i \leq l\}$. 
\end{lemma}
\begin{proof}
We argue by induction in the $n$-regular rank of $\Gamma$ that $S$ can be coded in the sorts $\Gamma \cup \{ \Gamma/\Delta_{i} \ | \ i \leq l\}$. % If $S$ has a minimum element $m$, this element is a code of such definable set. Thus we may assume that $S$ does not have a minimal element.\\
For the base case, we suppose that $\Gamma$ is $n$-regular. We first assume that $\Gamma$ is dense, and we aim to prove that $\beta$ and $\ulcorner S \urcorner$ are interdefinable. It is clear that $\ulcorner S \urcorner \in dcl^{eq}(\beta)$. For the converse let $\sigma$ be any automorphism of the monster model $\mathfrak{M}$ and suppose that $\sigma(\beta) \neq \beta$. Without loss of generality, $\beta < \sigma(\beta)$. By density we can find $n$-elements in the interval $(\beta, \sigma(\beta))$. By $n$-regularity and density there is an  element $\delta$ such that $\beta < n\delta< \sigma(\beta)$. Thus $\sigma(S) \subsetneq S$. \\
 We now assume that $\Gamma$ is discrete and let $1$ be its minimal element. There is a unique natural number $0 \leq i \leq n-1$ such that $\beta+i$ is $n$-divisible, because $\{ \beta, \beta+1, \dots, \beta+(n-1)\}$ is an interval with at least $n$-elements. Let $i_{0}$ be the index such that $\beta+i_{0}$ is $n$-divisible. Then $x \in S$ if and only if $nx \geq \beta+i_{0}$. Thus $\frac{\beta+i_{0}}{n}$ is the minimal element of $S$ and thereby it is a code for $S$. \\ 

We proceed to show the inductive step, and we consider the sequence $\{0\}=\Delta_{0} < \Delta_{1} < \dots < \Delta_{l+1}=\Gamma$ witnessing that $\Gamma$ has $n$-regular rank equal to  $l+1$. Let $\rho_{\Delta_{1}}: \Gamma \rightarrow \Gamma/\Delta_{1}$ be the canonical projection map, and note that  $\Gamma/\Delta_{1}$ is an ordered abelian group of $n$-regular rank $l$. First we suppose that $\rho_{\Delta_{1}}(\beta)$ is not $n$-divisible. Then $S$ is interdefinable with $S_{\Delta_{1}}=\{ \eta \in \Gamma/\Delta_{1} \ | \ n\eta > \rho_{\Delta_{1}}(\beta)\}$. By the induction hypothesis, such end-segment can be coded in the sorts $\Gamma/ \Delta_{1} \cup \{ (\Gamma/\Delta_{1})/ (\Delta_{i}/\Delta_{1})  \ | \ 2 \leq i \leq l\}$. As each of the sorts $(\Gamma/\Delta_{1})/ (\Delta_{i}/\Delta_{1})$ can be canonically identified with $\Gamma/\Delta_{i}$, the conclusion of the statement follows. \\
We consider the case where $\rho_{\Delta_{1}}(\beta)$ is $n$-divisible, i.e. there is some $\eta \in \Gamma$ such that $n \rho_{\Delta_{1}}(\eta)=\rho_{\Delta_{1}}(\beta)$. Note that $\rho_{\Delta_{1}}(\eta)= \min (S_{\Delta_{1}})$. If $\Delta_{1}$ is discrete, then $S$ has a minimum and this minimal element is a code for $S$. Thus without loss of generality  $\Delta_{1}$ is dense. We aim to show that $\beta$ and $\ulcorner S \urcorner$ are interdefinable. In fact, let $\sigma$ be an automorphism of the monster model $\mathfrak{M}$ fixing 
$\ulcorner S \urcorner$. We want to show that it fixes also $\beta$. We argue by contradiction, and we assume that $\beta<\sigma(\beta)$. As $\rho_{\Delta_{1}}(\beta)\in dcl^{eq}( \ulcorner S\urcorner)$, we have $\sigma(\beta)- \beta \in \Delta_{1}$. Fix some element $\eta \in \Gamma$ such that $n\eta +\Delta_{1}=\beta+\Delta_{1}$. We can find elements $\delta_{1}<\delta_{2} \in \Delta_{1}$ such that $\beta=n \eta+\delta_{1}$ and $\sigma(\beta)=n\eta+\delta_{2}$. By $n$-regularity and density of $\Delta_{1}$ we can find an element $\gamma \in \Delta_{1}$ such that  $\delta_{1}< n\gamma < \delta_{2}$, so we have $ \beta< n(\gamma+\eta)< \sigma(\beta)$ and hence $S \subsetneq \sigma(S)$, as desired. $\square$
\end{proof}
\begin{proposition} \label{codesegment}Let $\Gamma$ be an ordered abelian group of bounded regular rank, and let $S \subseteq \Gamma$ be a definable end-segment. Then $\ulcorner S \urcorner$ is interdefinable with a tuple of elements in the sorts $\Gamma \cup \{ \Gamma/ \Delta \ | \ \Delta \in RJ(\Gamma)\}$.  Consequently, any initial segment  is also coded in the sorts $\Gamma \cup \{ \Gamma/ \Delta \ | \ \Delta \in RJ(\Gamma)\}$.
\end{proposition}
\begin{proof}
By Proposition \ref{charend} it is sufficient to code divisibility end-segments. We may assume that $S= \{ \eta \in \Gamma \ | \ n \eta +\Delta \geq \beta+\Delta\}$. Therefore $S$ is interdefinable with $S_{\Delta}=\{ z \in \Gamma/\Delta \ | \ nz \geq \rho_{\Delta}(\beta)\}$; this is a definable end-segment of $\Gamma/\Delta$. The statement follows immediately from Lemma \ref{cociente} combined with Fact \ref{quotientsareok}. \\
The second part of the statement follows by noticing that any initial segment is the complement of an end-segment. $\square$
\end{proof}

\section{An abstract criterion to eliminate imaginaries}\label{criterion}
The following is \cite[Lemma 1.17]{criteria}.
\begin{theorem}\label{weakEI}Let $T$ be a first order theory with home sort $K$. Let $\mathcal{G}$ be some collection of sorts. If the following conditions all hold, then $T$ has weak elimination of imaginaries in the sorts $\mathcal{G}$. 
\begin{enumerate}
\item \emph{Density of definable types:} for every non-empty definable set $X \subseteq K$ there is an $acl^{eq}(\ulcorner X \urcorner )$-definable type in $X$.
\item \emph{Coding definable types:} every definable type in $K^{n}$ has a code in $\mathcal{G}$ (possibly infinite). That is, if $p$ is any (global) definable type in $K^{n}$, then the set $\ulcorner p \urcorner$ of codes of the definitions of $p$ is interdefinable with some (possibly infinite) tuple from $\mathcal{G}$.
\end{enumerate}
\end{theorem}
\begin{proof}
 A very detailed proof can be found in \cite[Theorem 6.3]{will}. The first part of the proof shows weak elimination of imaginaries as it is shown that for any imaginary element $e$ we can find a tuple $a \in \mathcal{G}$ such that $e \in dcl^{eq}(a)$ and $a \in acl^{eq}(e)$. $\square$ 
\end{proof}
We will use this criterion to prove that any pure ordered abelian group with bounded regular rank admits weak elimination of imaginaries once the quotient sorts are added. 
\subsection*{Coding of definable types}\label{codingdef}
In this subsection we show that any definable type $p(\mathbf{x})$ can be coded in the quotient sorts.
\begin{proposition}\label{codedef} Let $\Gamma$ be an ordered abelian group and $p(\mathbf{x}) \in S_{n}(\Gamma)$ be a definable type. Then $p(\mathbf{x})$ can be coded in the quotient sorts. 
\end{proposition}
\begin{proof}
Let $p(\mathbf{x})$ be a definable type in $n$ variables over $\Gamma$.  By quantifier elimination (Theorem \ref{QEbounded}), $p(\mathbf{x})$ is completely determined by formulas of the following forms:
\begin{itemize}
\item \emph{First kind:} 
\begin{align*}
\phi_{1}(\mathbf{x}, \beta):=\displaystyle{\sum_{i\leq n} z_{i} x_{i} +\Delta <   \beta+k +\Delta}
\end{align*}
or 
\begin{align*}
\psi_{1}(\mathbf{x}, \beta):=\displaystyle{\sum_{i\leq n} z_{i} x_{i} +\Delta > \beta+k+\Delta}
\end{align*}
where $\beta \in \Gamma$,  $\Delta \in RJ(\Gamma)$ and $k, z_{i} \in \mathbb{Z}$. 
\item  \emph{Second kind:}  
\begin{align*}
\phi_{2}(\mathbf{x}, \beta):=\displaystyle{\sum_{i\leq n} z_{i} x_{i}\equiv_{\Delta+ l\Gamma} \beta+k}
\end{align*}
where $\beta \in \Gamma$,  $\Delta \in RJ(\Gamma)$, $k, z_{i} \in \mathbb{Z}$ and $l \in \mathbb{N}_{\geq 2}$.
\item \emph{Third kind:}
\begin{align*}
\phi_{3}(\mathbf{x},\beta):=\displaystyle{\sum_{i\leq n} z_{i} x_{i}\equiv_{\Delta} \beta+k}
\end{align*}
where $\beta \in \Gamma$,  $\Delta \in RJ(\Gamma)$, and $z_{i} \in \mathbb{Z}$.
\end{itemize}
The set  $\{ \beta \in \Gamma \ | \ \phi_{1}(\mathbf{x},\beta) \in p(\mathbf{x}) \}$ is an end-segment of $\Gamma$, so it can be coded in the quotient sorts by Proposition \ref{charend} and \ref{codesegment}. Likewise, the set $\displaystyle{\{\beta \in \Gamma \ | \ \psi_{1}(\mathbf{x},\beta) \in p(\mathbf{x}) \}}$ is an initial segment of $\Gamma$, and it admits a code in the quotient sorts.\\
Let $X=\{ \beta \in \Gamma \ | \ \phi_{2}(\mathbf{x},\beta) \in p(\mathbf{x})\}$, then $X$ is either empty or we can take $\beta_{0} \in X$ and $\ulcorner X \urcorner$ is interdefinable with $\pi_{\Delta}^{l}(\beta_{0}) \in \Gamma/(\Delta+l\Gamma)$.\\
Lastly, $\displaystyle{Z=\{ \beta \in \Gamma \ | \ \phi_{3}(\mathbf{x},\beta) \in p(\mathbf{x})\}}$ is either empty or for any element $\beta_{0} \in Z$, we have that $\ulcorner Z \urcorner$ is interdefinable with $\rho_{\Delta}(\beta_{0}) \in \Gamma/ \Delta$. $\square$ \\
\end{proof}
\subsection*{Density of definable types}\label{densitydef}
In this subsection we prove the first condition required in Hrushovski's criterion: the density of definable types in algebraically closed sets.\\
The following will be a useful fact to obtain our result.
\begin{fact}\label{bueno}Let $X\subseteq \Gamma$ be a  definable set without a minimum element. Then there is a $\ulcorner X \urcorner$ definable end-segment $S$ such that $X$ is co-initial in $S$. 
\end{fact}
\begin{proof}
Let $I=\{ \beta \in \Gamma \ | \ (-\infty, \beta]\cap X=\emptyset\}$. $I$ is a $\ulcorner X \urcorner$-definable initial segment of $\Gamma$. Let $S=\Gamma \backslash I$, it is sufficient to verify that $X$ is co-initial in $S$. Let $\beta \in S$, then $(-\infty, \beta] \cap X \neq \emptyset$. Because $X$ does not have a minimum, we can find an element $x \in X$ such that $x < \beta$, as required. $\square$
\end{proof}

\begin{proposition}\label{density} Let $\Gamma$ be an ordered abelian group of bounded regular rank and $X \subseteq \Gamma$ a definable set. There is a global type $p(x) \vdash x \in  X$ such that $p(x)$ is definable over $\acl^{eq} (\ulcorner X \urcorner)$. 
\end{proposition}
\begin{proof}
Let $X \subseteq \Gamma$ be a $1$-definable set. If $X$ has a minimum element $a$, the statement follows immediately by taking the type of this element. Thus we may assume that $X$ does not have a minimum, by Fact \ref{bueno} there is a $\ulcorner X \urcorner$-definable end-segment $S$ such that $X$ is co-initial in $S$. In particular the type:
\begin{align*}
    \Sigma_{S}^{gen}(x)=\{ x \in S \cap X\} \cup \{ x \notin B \ | \ B\subsetneq S \ \text{and $B$ is a definable end-segment}\}
\end{align*}
is a consistent partial type which is $\ulcorner S \urcorner$-definable. \\
Let $\pi: \mathbb{N} \rightarrow \mathbb{N}\times \mathbb{N}_{\geq 1}$ be some fixed bijection. We now build by induction an increasing sequence of partial consistent types $(\Sigma_{i}(x) \ | \ i < \omega)$ in the following way:
\begin{itemize}
    \item \textbf{Stage $0$:} Set $\Sigma_{0}(x)=\Sigma_{S}^{gen}(x)$,
    \item \textbf{Stage $i+1$:} Let $\pi(i)=(k,l)$, at this stage we want to decide the congruence modulo the subgroup $\Delta_{k}+l\Gamma$. To keep the notation simple we assume that $l\geq 2$ and we use the projection map $\pi_{\Delta_{k}}^{l}:=\Gamma \rightarrow \Gamma/ (\Delta_{k}+l\Gamma)$. If $l=1$ we argue in the same manner to fix the coset of $\Delta_{k}$ and instead we use the projection map $\rho_{\Delta_{k}}:\Gamma \rightarrow \Gamma/\Delta_{k}$.\\
    We proceed by cases:
    \begin{itemize}
        \item [a)]Set $\Sigma_{i+1}(x)=\Sigma_{i}(x)\cup \{ \pi_{\Delta_{k}}^{l}(x) \neq \pi_{\delta_{k}}^{l}(\beta) \ | \ \beta \in \Gamma\}$ if it is consistent. 
        \item [b)] Otherwise, let $A_{i}=\{\eta_{1},\dots,\eta_{r_{i}}\} \subseteq \Gamma/(\Delta_{k}+l\Gamma)$ be the finite set of cosets such that $\Sigma_{i}(x)\cup\{ \pi_{\Delta_{k}}^{l}(x)=\eta_{j}\}$ is consistent. Take an element $\hat{\eta}\in A_{i}$ and set $\Sigma_{i+1}(x)=\Sigma_{i}(x) \cup \{ \pi_{\Delta_{k}}^{l}(x)=\hat{\eta}\}$.
    \end{itemize}
    Let $\mathfrak{M}$ be the monster model and
    \begin{align*}
    \mathcal{J}&=\{ i \in \mathbb{N} \ | \ \Sigma_{i}(x)\cup \{ \pi_{\Delta_{k}}^{l}(x) \neq \pi_{\delta_{k}}^{l}(\beta) \ | \ \beta \in \Gamma\} \ \text{is inconsistent}\}
    \end{align*}
\begin{claim}\label{lindo} For any $\sigma \in Aut(\mathfrak{M}/ \acl^{eq}(\ulcorner X \urcorner))$ the following conditions hold:
\begin{enumerate}
    \item For any $i \in \mathbb{N}$ $\sigma(\Sigma_{i}(x))=\Sigma_{i}(x)$ and
    \item For any $i \in \mathcal{J}$ $\sigma(A_{i})=A_{i}$.
\end{enumerate}
In particular, as $\sigma$ is arbitrary, then $A_{i}\subseteq \acl^{eq}(\ulcorner X \urcorner)$. 
\end{claim}
\begin{proof}
We argue by induction on $i$ to show that for any $\sigma \in Aut(\mathfrak{M}/\acl^{eq}(\ulcorner X \urcorner))$ we have that $\sigma(\Sigma_{i}(x))=\Sigma_{i}(x)$ and if $i \in \mathcal{J}$ then $\sigma(A_{i})=A_{i}$.\\
For the base case, fix some $\sigma \in Aut(\mathfrak{M}/\acl^{eq}(\ulcorner X \urcorner))$. Then $\sigma(\Sigma_{0}(x))=\Sigma_{0}(x)$ because $\Sigma_{S}^{gen}(x)$ is $\ulcorner S \urcorner$-definable and $\ulcorner S \urcorner \in \dcl^{eq}(\ulcorner X \urcorner)$.\\
Suppose the statement holds for $i$ and fix some $\sigma \in Aut(\mathfrak{M}/\acl^{eq}(\ulcorner X \urcorner)).$\\
If $\Sigma_{i+1}(x)=\Sigma_{i}(x)\cup \{ \pi_{\Delta_{k}}^{l}(x) \neq \pi_{\Delta_{k}}^{l}(\beta) \ | \ \beta \in \Gamma\}$, then
\begin{align*}
\sigma(\Sigma_{i+1}(x))&= \sigma(\Sigma_{i}(x)) \cup\{ \pi_{\Delta_{k}}^{l}(x)\neq \pi_{\Delta_{k}}^{l}(\sigma(\beta)) \ | \ \beta \in \Gamma\}\\
&=\Sigma_{i}(x)\cup\{ \pi_{\Delta_{k}}^{l}(x)\neq \pi_{\Delta_{k}}^{l}(\sigma(\beta)) \ | \ \beta \in \Gamma\}=\Sigma_{i+1}(x).
\end{align*}
Let's assume that $\Sigma_{i+1}(x)=\Sigma_{i}(x) \cup \{ \pi_{\Delta_{k}}^{l}(x)=\eta\}$ for some $\eta \in A_{i}$. We first argue that $\sigma(A_{i})=A_{i}$. By definition of $A_{i}$:\\
\begin{align*}
\mu \in A_{i}& \ \text{if and only if} \ \Sigma_{i}(x) \cup \{ \pi_{\Delta_{k}}^{l}(x)=\mu \} \ \text{is consistent.}\\
\end{align*}
Let $\mu \in A_{i}$, then $\displaystyle{\Sigma_{i}(x) \cup \{ \pi_{\Delta_{k}}^{l}(x)=\mu\}}$ is consistent. As $\sigma$ is an automorphism, then 
\begin{align*}
\sigma(\Sigma_{i}(x)) &\cup \{ \pi_{\Delta_{k}}^{l}(x)=\sigma(\mu)\} \ \text{is consistent.}\\
\end{align*}
By the induction hypothesis, $\sigma(\Sigma_{i}(x))=\Sigma_{i}(x)$. Hence:
\begin{align*}
\Sigma_{i}(x) \cup \{ \pi_{\Delta_{k}}^{l}(x)=\sigma(\mu)\}\ \ \text{is consistent.}\\
\end{align*}
Consequently, $\sigma(\mu) \in A_{i}$ and we conclude that $\sigma(A_{i}) \subseteq A_{i}$. We argue in a similar manner with $\sigma^{-1}$ to show that $A_{i}\subseteq \sigma(A_{i})$.\\
As for any $\sigma \in Aut(\mathfrak{M}/ \acl^{eq}(\ulcorner X \urcorner))$, $\sigma(A_{i})=A_{i}$ and $A_{i}$ is a finite set, then $ A_{i} \subseteq \acl^{eq}(\ulcorner X \urcorner)$. 
In particular, $\eta \in \acl^{eq}(\ulcorner X \urcorner)$ where $\Sigma_{i+1}(x)=\Sigma_{i}(x)\cup \{ \pi_{\Delta_{k}}^{l}(x)=\eta \}$.  Then for any $\sigma \in Aut(\mathfrak{M}/\acl^{eq}(\ulcorner X \urcorner))$ we have that $\displaystyle{\sigma(\Sigma_{i+1}(x))=\Sigma_{i+1}(x)}$, as required. 
\end{proof}
\end{itemize}
Let $\displaystyle{\Sigma_{\infty}(x)=\bigcup_{i \in \mathbb{N}} \Sigma_{i}(x)}$, this is a partial consistent type and $\Sigma_{\infty}(x) \vdash x \in X$. By quantifier elimination $\Sigma_{\infty}(x)$ determines a complete type $p(x)$. Then $p(x) \vdash x \in X$, and $p(x)$ is $acl^{eq}(\ulcorner X \urcorner)$-definable because $p(x)$ is completely determined by the data in $\Sigma_{\infty}(x)$, which is definable over $\acl^{eq}(\ulcorner X \urcorner)$ by Claim \ref{lindo} .$\square$
\end{proof}
\section{Main Results}\label{results}

\begin{theorem}\label{weakEIbounded} Let $\Gamma$ be an ordered abelian group of bounded regular rank (equivalently with finite spines). Then $\Gamma$ admits weak-elimination of imaginaries in the language $\mathcal{L}_{Q}$, once the quotient sorts are added.
\end{theorem}
\begin{proof}
By Theorem \ref{weakEI} it is sufficient to check that we have density of definable types and that we can code definable types in the quotient sorts. The first condition is Proposition \ref{density} and the second one is Proposition \ref{codedef}. $\square$
\end{proof}

\subsection*{The dp-minimal case}
In this section we show that a better statement can be achieved for the $dp$-minimal case. 
\begin{definition} Let $\Gamma$ be an ordered abelian group and $H$ some definable subgroup. A subset $\mathcal{C} \subseteq \Gamma$ is said to be a \emph{complete set of representatives modulo $H$} if:
\begin{enumerate}
\item given any $\gamma \in \Gamma$ there is some $\beta \in \mathcal{C}$ such that $\gamma-\beta \in H$. 
\item for any $\beta \neq \eta \in \mathcal{C}$ we have that  $\beta+H\neq \eta+H$. 
\end{enumerate}
\end{definition}
\begin{fact}\label{complete} Let $\Gamma$ be an ordered abelian group, $\Delta$ a convex subgroup and $k \in \mathbb{N}$. Let $\mathcal{C}$ be a complete set of representatives of $\Gamma$ modulo $k\Gamma$, then some subset $\mathcal{C}_{0} \subseteq \mathcal{C}$ is a complete set of representatives modulo $\Delta+k\Gamma$. 
\end{fact}
\begin{proof}
Let $\mathcal{C} \subseteq \Gamma$ be a complete set of representatives of $\Gamma$ modulo $k\Gamma$ and $\pi_{\Delta}^{k}: \Gamma \rightarrow \Gamma/(\Delta+k\Gamma)$ be the projection map. $\pi_{\Delta}^{k}(\mathcal{C})=\Gamma/(\Delta+k\Gamma)$, because for any $\gamma \in \Gamma$, there is some $\beta \in \mathcal{C}$ such that $\gamma-\beta \in k\Gamma$, in particular $\gamma-\beta \in \Delta+k\Gamma$. For each coset $\eta \in \Gamma/(\Delta+k\Gamma)$ choose an element $c_{\eta}\in \mathcal{C}$ such that $\pi_{\Delta}^{k}(c_{\eta})=\eta$. The set $\mathcal{C}_{0}=\{ c_{\eta}\ | \ \eta \in \Gamma/(\Delta+k\Gamma)\}$ is a complete set of representatives. 
\end{proof}
By Proposition \ref{dpmin}, an ordered abelian group is $dp$-minimal if and only if it does not have singular primes, i.e. for any $p$ prime number $[\Gamma:p\Gamma]< \infty$.  We consider the language $\mathcal{L}_{dp}$ extending $\mathcal{L}_{Q}$, where for each $k \in \mathbb{N}_{\geq 2}$ we add constants for the elements of the finite groups $\Gamma/k\Gamma$.

\begin{corollary} Let $\Gamma$ be a dp-minimal ordered group. Then $\Gamma$ admits elimination of imaginaries in the language $\mathcal{L}_{dp}$, where the quotient sorts are added. 
\end{corollary}
\begin{proof}
By Theorem \ref{weakEIbounded} and Fact \ref{all} it is sufficient to show that we can also code finite sets. Let $\Delta$ definable convex subgroup and $k \in \mathbb{N}$, the group $\Gamma/(\Delta+k\Gamma)$ is also finite. We first argue that $\Gamma/(\Delta+k\Gamma) \subseteq \dcl(\emptyset)$. Consider the $\emptyset$-definable function
\begin{align*}
    f:& \Gamma/k\Gamma \rightarrow \Gamma/(\Delta+k\Gamma)\\
    & \gamma+k\Gamma \rightarrow \gamma+(\Delta+k\Gamma).
\end{align*}
By Fact \ref{complete} $f$ is surjective.\\
Hence, it is enough to prove that finite sets of tuples in $\mathcal{S}=\{ \Gamma/\Delta \ | \ \Delta \in RJ(\Gamma)\}$ can be coded in the quotient sorts. As each of the sorts $\Gamma/\Delta$ is linearly ordered, there is a definable order induced over the finite products of quotients of $\Gamma/\Delta$, and thereby any finite set of tuples in $\mathcal{S}$ is already coded in $\mathcal{S}$. $\square$
\end{proof}

\newpage
\bibliographystyle{plain}

\end{document}